\documentclass[12pt,letterpaper]{amsart}
\usepackage{graphicx}
\usepackage{amsmath, amssymb, amsthm, mathrsfs}

\usepackage{setspace,verbatim}
\newtheorem{lemma}{Lemma}

\newcommand{\R}{\ensuremath{\mathbb{R}} }

\newcommand{\Di}{\ensuremath{\Delta^\infty} }
\newcommand{\Dio}{\ensuremath{\Delta^\infty_0} }

\newtheorem{theorem}{Theorem}

\newtheorem{prop}{Proposition}

\newcommand{\Dho}{\ensuremath{\Delta^H_0} }
\newtheorem*{sperner}{Sperner's Lemma}

\newtheorem*{sch}{Schauder's Theorem}

\newcommand{\ignore}[1]{}

\begin{document}

\title[A Fixed Point Theorem 
for the Infinite-dimensional Simplex]{
A Fixed Point Theorem\\
for the Infinite-dimensional Simplex}

\author[Douglas Rizzolo]{Douglas Rizzolo$^*$}
 
\author[Francis Edward Su]{Francis Edward Su$^{**}$}

\thanks{$^{*}$Harvey Mudd College, Claremont, CA  91711.
E-mail:{\tt drizzolo@hmc.edu}}

\thanks{$^{**}$Corresponding author. 
Department of Mathematics, Harvey Mudd College, Claremont, CA  91711.
Voice: 909-607-3616, Fax: 909-621-8366, E-mail: {\tt
  su@math.hmc.edu}.}

\thanks{The authors gratefully acknowlege 
partial support by NSF Grant DMS-0301129 (Su), as well as helpful
conversations with Jon Jacobsen.}

\dedicatory{Manuscript, Oct 23, 2006}

\subjclass[2000]{Primary 54H25; Secondary 47H10, 55M20}
\keywords{Schauder fixed point theorem, Brouwer fixed point theorem,
  Sperner's lemma, infinite-dimensional simplex}

\begin{abstract}
We define the infinite dimensional simplex to be the closure of the
convex hull of the standard basis vectors in $\R^\infty$, and prove
that this space has the {\em fixed point property}: any continuous
function from the space into itself has a fixed point.  Our proof is
constructive, in the sense that it can be used to find an approximate
fixed point; the proof relies on elementary analysis and Sperner's
lemma.  The fixed point theorem is shown to imply Schauder's fixed
point theorem on infinite-dimensional compact convex subsets of 
normed spaces.
\end{abstract}

\maketitle
                                                                               
\section{Introduction}

In finite dimensions, one of the simplest methods for proving the
Brouwer fixed point theorem is via a combinatorial result known as
Sperner's lemma \cite{Sper28}, which is a statement about labelled
triangulations of a simplex in $\R^n$.  In this paper, we use
Sperner's lemma to prove a fixed point theorem on an 
infinite-dimensional simplex in $\R^\infty$.  We also show that this
theorem implies the infinite-dimensional case of Schauder's fixed point
theorem on normed spaces.

Since $\R^\infty$ is locally convex, our theorem 
is a consequence of Tychonoff's fixed point theorem \cite{Smar74}.
However, some notable advantages of our approach are: 
(1) the constructive nature of Sperner's lemma provides a method for
producing approximate fixed points for functions on the
infinite-dimensional simplex, (2) the proof is based on elementary
methods in topology and analysis, and (3) our proof provides another
route to Schauder's theorem.

Fixed point theorems and their constructive proofs 
have found many important applications, ranging from proofs 
of the Inverse Function Theorem \cite{Lang97}, to proofs of the
existence of equilibria in economics \cite{Todd76, Yang99}, to the 
existence of solutions of differential equations \cite{Brow93, Smar74}.

\section{Working in $\R^\infty$}
Let $\R^\infty$ and $I^\infty = \prod [0,1]$ be the product of 
countably many copies of $\R$, and $I=[0,1]$, respectively.  
We equip $\R^\infty$ with the standard product topology, which is 
metrizable \cite{BePe75} by the complete metric
\[
\bar d(x,y) = \sum_{i=1}^\infty \frac{|x_i - y_i|}{2^i(1+|x_i-y_i|)}.
\]
In $\R^n$, a $k$-dimensional simplex, or \textit{$k$-simplex}, $\sigma^k$
is the convex hull of $k+1$ affinely independent points.
The \textit{standard $n$-simplex} in $\R^{n+1}$, denoted $\Delta^n$,
is the convex hull of the $n+1$ standard basis vectors of $\R^n$. 

The natural extension of this definition to $\R^\infty$ is 
to consider $\Delta^\infty$, the convex hull of the standard basis
vectors $\{e_i\}$ in $\R^\infty$, where $(e_i)_j= \delta_{ij}$, the
Kronecker delta function.  As convex combinations are finite
sums, this convex hull is:
\[
\Di = \{x\in \R^\infty | \sum_{i=1}^\infty x_i = 1,\  0\leq x_i \leq
1 \mbox{, and only finitely many $x_i$ are non-zero}\}.
\]
Unfortunately, $\Di$ is not closed; under 
the metric $\bar d$ the sequence $\{e_i\}$ converges to $\mathbf{0}$, 
which is not in $\Di$.
So consider, instead $\Dio$, the closure of $\Di$, which can be shown to be:
\[  
\Dio  = \{x\in \R^\infty | \sum_{i=1}^\infty x_i \leq 1 \mbox{ and }
0\leq x_i\leq 1\}.
\]
It is easy to see that $\Dio$ is convex.  It is also the closure
of the convex hull of 
the standard basis vectors $\{e_i\}$ and $\mathbf{0}$.
   It is also compact because it is a closed subset of $I^\infty$, 
   which is compact by Tychonoff's Theorem.
      \footnote{If one would like to avoid the Axiom
	of Choice, which is equivalent to Tychonoff's Theorem, it is not
	difficult to show that $I^\infty$ is a closed and totally bounded
	subset of the complete space $\R^\infty$, which implies
	compactness.}  
We call $\Dio$ the {\em standard infinite-dimensional simplex}.

It will be important for our purposes later to consider $F^n$, 
   the $n$-dimensional face of $\Dio$ given by $F^n = conv\{e_1,e_2,\dots,
	e_{n+1}\}$. Notice that each $F^n$ is closed and thus compact.

\section{Some Preliminary Machinery}

Let $\sigma^k = conv(x_0,...,x_k)$ be a $k$-simplex in $\R^n$.   
Let $T$ be a triangulation of $\sigma^k$ and $V$ be the set of
vertices of $T$ (i.e., the vertices of simplices in $T$).   
A \emph{Sperner labelling} of the triangulation $T$ 
is a labelling function $\ell:V\rightarrow \{0,\dots,k\}$ such that 
\[ 
\mbox{ if } J \subseteq \{0,\dots,k\} \mbox{ and } v \in conv\{x_j | j \in J \}
\mbox{, then } h(v) \in J.
\]
A $k$-simplex $\tau$ of $T$ 
is called a {\em fully-labelled simplex} (or {\em full})
if the image of the vertices of $\tau$ under $\ell$ maps onto
$\{0,\dots,k\}$.  Note that $\tau$ has exactly $k+1$ vertices, so all
the vertices have distinct labels.

\begin{sperner} 
Let $\sigma^k$ be a $k$-simplex in $\R^n$ with triangulation $T$ and let
$\ell$ be a Sperner-labelling of $T$.
Then the number of full simplices of $T$ is odd (and hence, non-zero). 
\end{sperner}

Though we will not prove this theorem here, an exposition of such 
proofs can be found in
\cite{Su99}.  In particular, there are constructive
``path-following'' proofs that locate the full simplex by tracing a
path of simplices through the triangulation.  Such path-following
proofs have formed the basis of
algorithms for locating fixed points of functions in 
finite-dimensional spaces, e.g., see \cite{Todd76} for a nice survey.
In Section \ref{sec:fixed-dio}, we show how to use Sperner's lemma for
a fixed point theorem in the infinite-dimensional space $\Dio$.

Another crucial theorem for our purposes states that, under
appropriate hypotheses, the existence of approximate fixed points
implies the existence of fixed points.  On the metric space $(X,d)$,
we can quantify the notion of an approximate fixed point by  
defining an \textit{$\epsilon$-fixed point}, which for a given function
$f$ is a point $x\in X$ such that $d(x,f(x))<\epsilon$.   
Versions of the following lemma may be found in, e.g., \cite{DuGr82, Smar74}.

\begin{lemma}[Epsilon Fixed Point Theorem] 
\label{le epsilonfixed}
 Suppose that $A$ is a
  compact subset of the metric space $(X,d)$ and that $f:A\rightarrow
  A$ is continuous.   If $f$ has an $\epsilon$-fixed point for every
  $\epsilon > 0$ then $f$ has a fixed point. 
\end{lemma}

\begin{proof} Let $\{a_n\}$ be a sequence of $1/n$-fixed points.
  That is, $d(a_n,f(a_n)) < 1/n$ for all $n$.   Since $A$ is compact
  it is sequentially compact and thus $\{a_n\}$ has a convergent
  subsequence, which we denote $\{a'_n\}$ with $a_n'\rightarrow x \in
  A$.   Let $\epsilon >0$.   Since $a_n'\rightarrow x$ there exists
  $N_1$ such that $n\geq N_1$ implies that $d(a'_n,x) < \epsilon /2$.
  Let $N = \max ( N_1, 2/\epsilon)$.   Then $n\geq N$ implies that
\[
d(x, f(a'_n)) \leq d(x,a'_n) + d(a'_n,f(a'_n)) < \epsilon,
\]
so that $f(a'_n) \rightarrow x$.   However, since $f$ is continuous,
we also know that $f(a'_n)\rightarrow f(x)$.   Since limits are
unique, we conclude that $f(x)=x$, which completes the
proof. 
\end{proof}

Later it will be desirable to have an isometry between
$\Delta^{n-1}$, the standard $(n-1)$-simplex in $\R^{n}$, and
$F^{n-1}$.  The easiest way to do this is to consider
$\R^n$ as a subspace of $\R^\infty$ by projection onto the first $n$
factors, and restricting the metric on $\R^\infty$ to $\R^{n}$.  Call
this metric $\bar d_{n}$ and consider $\Delta^{n-1}$ 
in the metric space $(\R^{n}, \bar d_{n})$.
It is worthwhile to ensure that 
$(\R^{n}, \bar d_n)$ has a rich supply of continuous functions.  
Before proceeding, recall that all norms on $\R^n$ are equivalent 
and thus essentially interchangeable; we now prove that $\bar d_n$ is
interchangeable with norm-induced metrics on bounded sets.

\begin{lemma} 
\label{le metricequivalence}
Let $A$ be a bounded subset of the normed space $(\R^{n}
  , \| \cdot \|_\infty)$.   On $A$, the metric $\bar d_n$ is
  equivalent to the metric induced by the norm $\| \cdot \|_\infty$.
\end{lemma}

\begin{proof} Suppose that $x,y \in \R^n$.   We see that 
\[
\bar d_n(x,y) = \sum_{i=1}^n \frac{|x_i - y_i|}{2^i(1+|x_i-y_i|)} \leq
n\|x-y\|_\infty.
\]
Now, since $A$ is bounded, there
is some $M$ such $\|x-y\|_\infty \leq M$ for $x,y \in A$.   Thus we
see that
\[\frac{\|x-y\|_\infty}{2^n(1+M)} \leq
\frac{\|x-y\|_\infty}{2^n(1+\|x-y\|_\infty)} \leq \bar d_n(x,y),
\]
which implies that 
\begin{equation}
  \label{eq:norm-bound}
\|x-y\|_\infty \leq 2^n(1+M)\bar d_n(x,y).
\end{equation}
Thus
$\bar d_n$ is equivalent to the metric induced by the norm on
$A$.\end{proof} 

Lemma \ref{le metricequivalence} tells us that bounded subsets of $\R^n$ have
the same continuous functions regardless of whether they are
considered as subsets of a normed space or as subsets of $(\R^n , \bar
d_n)$.  Importantly, notice that $\Delta^{n-1}$ is bounded.
Furthermore, the isometry 
$f:\Delta^{n-1} \rightarrow F^{n-1}$
between $\Delta^{n-1}$ in $(\R^n,\bar d_n)$
and $F^{n-1}$ in $\R^\infty$ is clearly given
by $f(x) = f(x_1,x_2,\dots,x_n) =
(x_1,x_2,\dots,x_n,0,0,\dots)$.  
This is important because it implies that $F^{n-1}$ has an arbitrarily
small barycentric subdivision.
Recall that the diameter of a set $X$ is
$d(X) = \sup_{x,y\in X} d(x,y)$ and if $\mathscr T$ is a family of
sets, then $size(\mathscr T) = \sup_{\sigma \in \mathscr T}
d(\sigma)$.
Thus, given $\epsilon >0$, $F^{n-1}$ has a barycentric subdivision $\mathscr
T$ with $size(\mathscr T) < \epsilon$. 

************************

Now we are ready to prove a fixed point theorem for $\Dio$.

\section{A Fixed Point Theorem for $\Dio$}
\label{sec:fixed-dio}

\begin{theorem}
\label{fixed-dio}
Suppose that $f:\Dio\rightarrow \Dio$ is continuous.   Then $f$ has a
fixed point. 
\end{theorem} 

\begin{proof} 
Since $\Dio$ is compact, by Lemma \ref{le epsilonfixed}, it is sufficient to show that $f$
has an $\epsilon$-fixed point for each $\epsilon >0$.   Let $\epsilon
>0$ be given.   Choose $N \geq \log_2(2/\epsilon)+1$.   Notice that for
$x,y \in \Dio$, this implies that 
\begin{equation}
\label{eq:eps-over-2-N+1on}
\sum_{i=N+1}^\infty   \frac{|x_i - y_i|}{2^i(1+|x_i-y_i|)}  \leq
\sum_{i=N+1}^\infty \frac{1}{2^i} < \frac{\epsilon}{2}.
\end{equation}
Since $f$ maps between countably infinite-dimensional spaces, 
we can write $f$ in terms of its components: 
$f(x) = (f_1(x),f_2(x),\dots)$.   
Since $f$ is continuous, 
$f_i$ is continuous for each $i$.   Consider the function 
\[
g(x) =(g_1(x),g_2(x),\dots) =  (f_1(x),f_2(x),\dots, f_N(x) ,
1-\sum_{i=1}^N f_i(x), 0,0,0,\dots ).
\]   
Since each $f_i$ is continuous and finite sums of continuous function
are continuous, $g_i$ is continuous for each $i$.   Furthermore, we
see that $g:F^N \rightarrow F^N$.   Consequently, $g$ is continuous.    

Let $\epsilon_0 = \frac{\epsilon}{8(N+1)}$ and $\epsilon_1 =
\frac{\epsilon}{2^{N+5}(N+1)}$.   Since $g$ is continuous on a compact
set, it is uniformly continuous.   Thus there exists $\delta_1 >0$
such that $\bar d(x,y) < \delta_1$ implies that $\bar d(g(x),g(y)) <
\epsilon_1$.   Let $\delta = \min ( \delta_1 , \epsilon_1)$.   Since
$F^N$ can be triangulated with an arbitrarily small triangulation, let
$\mathscr T$ be a triangulation with $size(\mathscr T) < \delta$.   
Label the vertices of $\mathscr T$ with the map 
\[
\ell(x) = \mbox{argmax}_i (x_i - g_i(x)).
\]
Recall that the \emph{argmax} function returns the index of the
largest element of the argument, and if there are multiple 
indices that give the maximum value, 
the argmax function returns the least of these indices.   

Observe that $\ell(x)$
produces a Sperner labeling on the vertices of $\mathscr T$.   Thus by
Sperner's Lemma, there exists a fully-labeled simplex in $\mathscr T$.   
This simplex can be found using the path-following method described in
\cite{Su99}.   Let $\{x^1,x^2, \dots x^{N+1}\}$ be the
vertices of this simplex where the index of each vertex is its Sperner
label.   From this, we see that for all $j$,
\[
x^i_i - g_i(x^i) \geq x^i_j - g_j(x^i). 
\]
Furthermore, since for each $x$ in $F^N$, we have
\[
\sum_{j=1}^{N+1} x_j = \sum_{j=1}^{N+1} g_j(x) =1,
\]
there is at least one $j$ such that $g_j(x) \leq x_j$.  
In particular, since $\ell(x^i)=i$, this implies that for each $x^i$,
\[ 
x^i_i - g_i(x^i) = \max_j ( x^i_j - g_j(x^i)) \geq 0.
\]
Since $size(\mathscr T) < \delta$ we have that, for all $i$, $\bar
d(x^1, x^i) < \delta$.   From the bound (\ref{eq:norm-bound}) in
Lemma \ref{le metricequivalence} (note in this case $M=1$ and $n=N+1$), 
we find that for all $i, j$,
\begin{equation}
\label{eq:delta-bound}
|x^1_j - x^i_j| < 2^{N+2}\delta \leq 2^{N+2} \epsilon_1 \leq \epsilon_0.
\end{equation}
By the same logic, we have that for all $i, j$,
\begin{equation}
\label{eq:delta-eps-bound}
|g_j(x^1) - g_j(x^i)| < 2^{N+2}\epsilon_1 \leq \epsilon_0.
\end{equation}
Consequently, we have that
\[
x^1_j + \epsilon_0 > x^i_j \quad \mbox{ and } \quad -g_j(x^i) <
\epsilon_0 - g_j(x^1)
\]
which, in turn, implies that 
\[
2\epsilon_0 + x^1_j - g_j(x^1) > x^i_j - g_j(x^i)
\]
for all $i$ and $j$.
In particular, this implies that the following list of inequalities
hold (simply let $i=j$ and run through all $i$): 
\[
  \begin{array}{cccc} 
   2\epsilon_0 + x^1_1 - g_1(x^1) & > &x^1_1 - g_1(x^1)& \geq 0 , \\
   2\epsilon_0 + x^1_2 - g_2(x^1) & > &x^2_2 - g_2(x^2)& \geq 0 , \\
  \vdots & & \vdots \\
  2\epsilon_0 + x^1_{N+1} - g_{N+1}(x^1) & > & x^{N+1}_{N+1} -
   g_{N+1}(x^{N+1})& \geq 0 .
\end{array} 
\]
Summing down each column yields the following inequality.
\[2\epsilon_0(N+1) + \sum_{i=1}^{N+1} x^1_i - \sum_{i=1}^{N+1}
g_i(x^1) > \sum_{i=1}^{N+1}\left( x_i^i-g_i(x^i)\right) \geq 0. 
\]
Now we recall that for all $i$, $x_i^i-g_i(x^i) \geq 0$ and 
\[
\sum_{i=1}^{N+1} x^1_i - \sum_{i=1}^{N+1} g_i(x^1)=1-1 =0.
\]
Consequently, 
\[
\begin{split} 
  2\epsilon_0(N+1)& = 2\epsilon_0(N+1)+ \sum_{i=1}^{N+1} x^1_i -
  \sum_{i=1}^{N+1} g_i(x^1) \\ 
  & >\sum_{i=1}^{N+1}\left( x_i^i-g_i(x^i)\right)\\
  &=\sum_{i=1}^{N+1}\left| x_i^i-g_i(x^i)\right|. 
\end{split}
\] 
Using (\ref{eq:delta-bound}) and (\ref{eq:delta-eps-bound}) 
and the continuity of $g$, for all $i$, we have that:
$
|x^1_i - g_i(x^1)| 
\leq |x^1_i - x^i_i| + |x^i_i - g_i(x^i)| + |g_i(x^i)-g_i(x^1)|
< 2 \epsilon_0 +  |x^i_i - g_i(x^i)|
$.
Hence,
\[ \begin{split} \bar d(x^1 , g(x^1)) = \sum_{i=1}^{N+1} \frac{|x^1_i - g_i(x^1)|}{2^i(1+|x^1_i - g_i(x^1)|)}  &  \leq \sum_{i=1}^{N+1} |x^1_i - g_i(x^1)| \\
 & <  \sum_{i=1}^{N+1} \left(2 \epsilon_0 +  |x^i_i - g_i(x^i)|\right) \\
 & < 4(N+1)\epsilon_0 \\
 & = \frac{\epsilon}{2}. \end{split} \]
Let $y = (x^1_1,x^1_2,\dots,x^1_N,0,0,0,\dots)$.   
We see that 
\begin{equation}
\label{eq:eps-over-2-1toN}
\begin{split}
  \sum_{i=1}^{N} \frac{|y_i - f_i(y)|}{2^i(1+|y_i - f_i(y)|)}
&=\sum_{i=1}^{N} \frac{|y_i - g_i(y)|}{2^i(1+|y_i - g_i(y)|)}\\ 
& = \sum_{i=1}^{N} \frac{|x^1_i - g_i(x^1)|}{2^i(1+|x^1_i - g_i(x^1)|)} \\
& \leq \sum_{i=1}^{N+1} \frac{|x^1_i - g_i(x^1)|}{2^i(1+|x^1_i - g_i(x^1)|)} \\
& < \frac{\epsilon}{2}.
\end{split}
\end{equation}
From (\ref{eq:eps-over-2-N+1on}) and (\ref{eq:eps-over-2-1toN}),
we have 
\[
  \begin{split} \bar d(y,f(y)) & = \sum_{i=1}^\infty 
      \frac{|y_i - f_i(y)|}{2^i(1+|y_i - f_i(y)|)} \\
  & =\sum_{i=1}^{N} \frac{|y_i - f_i(y)|}{2^i(1+|y_i - f_i(y)|)}
  +\sum_{i=N+1}^\infty   \frac{|y_i - f_i(y)|}{2^i(1+|y_i - f_i(y)|)} \\
  & < \frac{\epsilon}{2}+\frac{\epsilon}{2}\\
  &=\epsilon . 
\end{split}
\]
Therefore, $y$ is the desired $\epsilon$-fixed point.\end{proof}
Notice that the construction of the $\epsilon/2$ fixed point in $F^N$
in the proof above is identical to the construction of an $\epsilon/2$
fixed point for an arbitrary continuous function on $\Delta^N$,
because of the isometry between the two sets.   This construction, in
conjunction with Lemmas \ref{le epsilonfixed} and \ref{le
 metricequivalence}, provides a proof of the 
Brouwer Fixed Point Theorem on the finite-dimensional simplex, which is
similar to constructions found in, e.g., \cite{Todd76}.

\section{Schauder's Theorem}

A well-known infinite-dimensional fixed point theorem that holds for
normed spaces is Schauder's theorem \cite{DuGr82, Smar74}:

\begin{sch} Suppose that $X$ is a compact convex subset of the normed
  space $G$.   If $f:X\rightarrow X$ is continuous, then $f$ has a
  fixed point. \end{sch}

In this section we show how our proof of Theorem \ref{fixed-dio} can
be used to prove Schauder's Theorem for the case where $X$ is
infinite-dimensional.  (The finite-dimensional version of Schauder's
Theorem reduces to the Brouwer Fixed Point Theorem.)

Recall that a space $X$ has the \textit{fixed point property}.
if every continuous function $f:X\rightarrow X$ has a fixed point.
Note that this is a topological property, so if $X$ is homeomorphic to
$Y$ then $Y$ also has the fixed point property.
We will establish Schauder's theorem by noting that $\Dio$ is homeomorphic
to any infinite-dimensional compact convex subset of a normed space.

Define the vector space $H$ to be
\[H = \{x \in \R^\infty | \sum_{i=1}^\infty \frac{|x_i|}{2^i} < \infty \}.\]
It is not difficult to see that $H$ is indeed a vector space.
Furthermore, we see that $\|x\| = \sum_{i=1}^\infty \frac{|x_i|}{2^i}$
defines a norm on this space and the closure of the standard simplex in $H$ is 
\[
\Dho = \{x\in H | \sum_{i=1}^\infty x_i \leq 1 \mbox{ and } 0\leq
x_i\leq 1\} .
\]

\begin{prop} 
$\Dio$ is homeomorphic to $\Dho$. 
\end{prop}
The proof of this lemma is trivial using the homeomorphism $g:\Dio
\rightarrow \Dho$ being $g(x)=x$.
Note that $\Dho$ is an infinite-dimensional compact convex subset of
a normed space $H$.  Now consider the following proposition \cite{Klee55}:

\begin{prop} Every infinite-dimensional compact convex subset of a
  normed space is homeomorphic to the Hilbert Cube. 
\end{prop}

The significance of these propositions is that \emph{every}
infinite-dimensional compact convex subset of a normed space is
homeomorphic to $\Dio$.
Thus Theorem \ref{fixed-dio} implies the infinite-dimensional case of 
Schauder's Theorem.

\bibliographystyle{plain}	
\bibliography{topology}
\end{document}